\documentclass[12pt,reqno]{amsart}
\usepackage{amsmath,amssymb,amstext,amsfonts,epsfig,graphicx}
\usepackage[margin=1in]{geometry}
\usepackage{float}
\usepackage{graphicx}
\usepackage{subcaption}
\usepackage{epsfig}
\usepackage{amsmath}
\usepackage{amsfonts}
\usepackage{amssymb}
\usepackage{amsthm}

\title{Two Universality Properties Associated with the Monkey Model of Zipf's Law}
\author{Richard Perline}
\address{Flexible Logic Analytics, 34-50 80th Street, Jackson Heights, NY 11372, USA}
\email{richperline@gmail.com}
\author{Ron Perline}
\address{Department of Mathematics, Drexel University, Philadelphia, PA, 19104,USA}
\email{ronald.k.perline@drexel.edu}

\begin{document}

\begin{center}
{\tt A version of this paper was submitted to Entropy on 11/20/15.\\ This document was \TeX 'd on 11/27/15.}
\end{center}
\vskip .5in

\begin{abstract} 
The distribution of word probabilities in the monkey model of Zipf's law is associated with two universality properties:  
(1) the power law exponent converges strongly to $-1$
as the alphabet size increases and  the letter probabilities are specified as the spacings from a random division of the unit interval
for any distribution with a bounded 
density function on $[0,1]$; and (2), on a logarithmic scale the version of the model with a finite word length cutoff and unequal letter probabilities 
is approximately normally distributed in the part of the distribution away from the tails.
The first property is proved using a  remarkably general limit theorem  for the logarithm of sample spacings from Shao and Hahn,
and the second property follows from Anscombe's central limit theorem for a random number of i.i.d. random variables. 
The finite word length model leads to a hybrid Zipf-lognormal mixture distribution closely related to work in other areas.
\end{abstract}
\maketitle

\section{Introduction}

In his popular expository article on universality examples in mathematics and physics Tao 
\cite{Tao}
mentions Zipf's empirical power law of word frequencies
as lacking any ``convincing explanation for how the law comes about and why it is universal.''   By {\it universality} he means the idea where systems follow
certain macroscopic laws that are largely independent of their microscopic details.  In this article we drop down many levels from natural language
to discuss {\it two} universality properties associated with the toy monkey model of Zipf's law.  Both of these properties were considered in Perline \cite{Perline1996},
but the presentation there was incompletely developed, and  here we  clarify and greatly expand upon these two ideas by showing:  (1) how the model displays a nearly universal tendency towards a $-1$ exponent in its power law behavior; and (2) the significance of the central limit theorem (CLT) for a random number of i.i.d. variables in the model with a finite word length cutoff.  We will sometimes refer to the case with a finite word length cutoff as {\it Monkey Twitter}, as explained in Section 3. 

Our paper is laid out as follows. 
In Section 2 we prove the  strong tendency towards an approximate $-1$ exponent under very broad conditions by means of a limit theorem for the logarithms of random
spacing due to Shao and Hahn \cite{S&H}.
A somewhat longer approach to this proof is  in  Perline  \cite{Perline2015}, where an elementary derivation of power law behavior of the monkey  model is first given
before applying the Shao-Hahn limit result.
In Section 3, we  explain the underlying lognormal structure of the {\it central part} of the distribution of
word probabilities in the case of a {\it finite word length}  cutoff and show how it leads to a hybrid Zipf-lognormal distribution (called a lognormal-Pareto distribution
 in   \cite{Perline1996}). 
These are universality properties in exactly Tao's sense because the distribution of the word probabilities (the macro behavior) is within very broad bounds independent of the details of how the letter probabilities for the monkey's typewriter  are selected (the micro behavior).
In Section 4, we discuss the multiple connections between these results and  well-known research from other areas  where Pareto-Zipf type distributions
have been investigated.  
In the remainder of this section we sketch some historical background.     
 
It has been known for many years  that the monkey-at-the-typewriter scheme for generating random ``words'' conforms to an inverse power law of word frequencies
and therefore mimicks the inverse power form of Zipf's \cite{Zipf1949}  statistical regularity  for many natural languages. 
{\it However, Zipf's empirical word frequency law not only follows a power distribution, but  as he emphasized, it also very frequently exhibits 
an exponent in the neighborhood of $-1$.}   That is, letting $f_1\ge f_2\ge \dots f_r\dots\ge ...$ represent the observed ranked word frequencies, Zipf found $f_r\approx Cr^{-1}$,  with $C>0$ a constant depending on the sample size. 
Many qualifications have been raised regarding this approximation (including the issue of the divergence of the harmonic series); yet
as I.J. Good \cite{Good1969}  
 commented:   ``The Zipf law is unreliable but it is often a good enough approximation to demand an explanation.''   Figure 1 illustrates
the fascinating universal character of this word frequency law using the text from four different authors writing in four different European languages in four different centuries.
The plots are shown in the log-log form that Zipf employed, and they resemble countless other examples that researchers have studied since his time.
The raw text files for the four books used in the graphs were downloaded from the Public Gutenberg databank \cite{Gutenberg}.  

\begin{figure}[h]
\centering
\includegraphics[width=4.in]{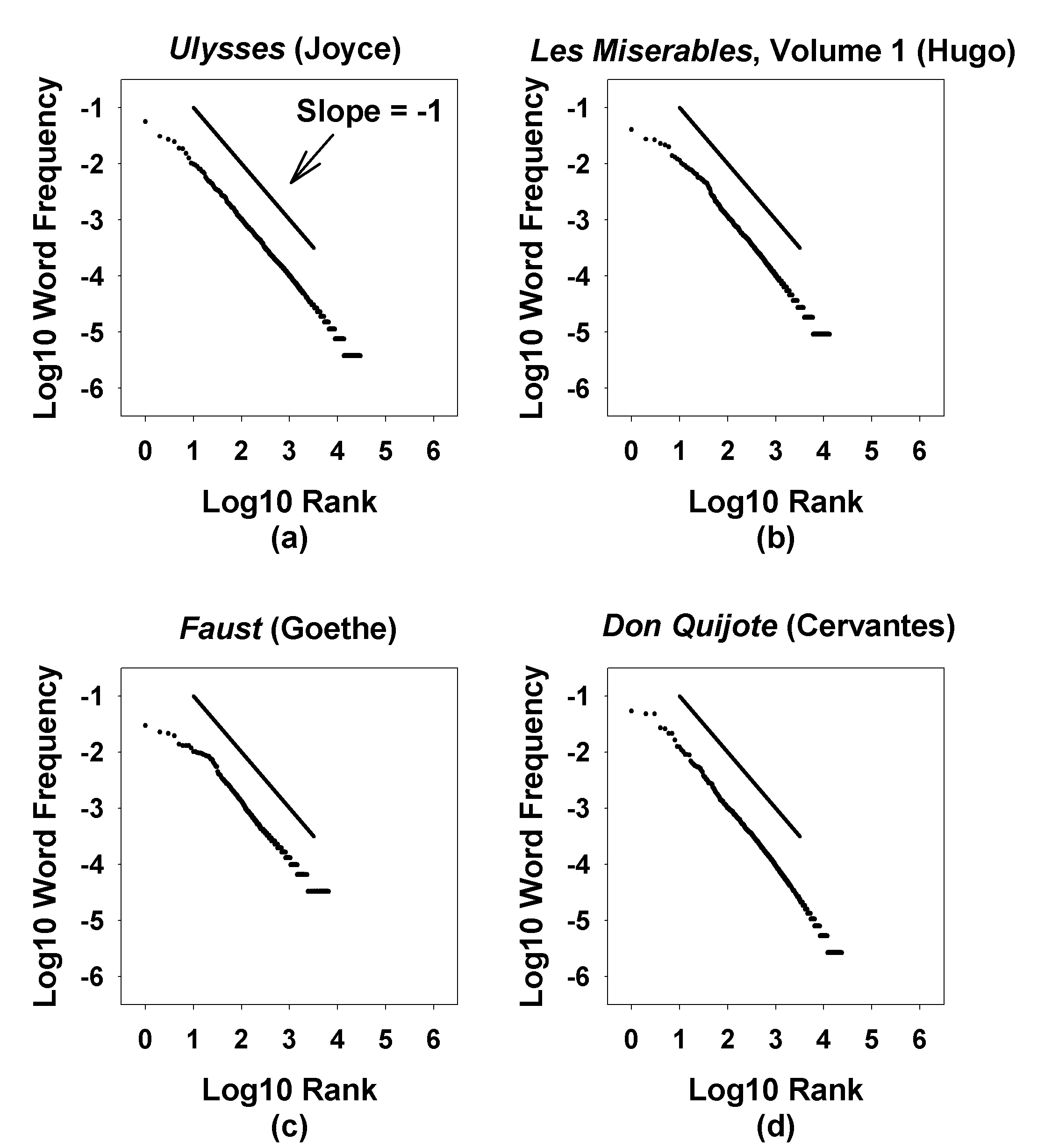}
\caption{Log-log plots of relative word frequencies by rank for four authors writing in four different European languages in four different centuries.  The 
approximate -1 slope  in all the graphs is an iconic feature of Zipf's word frequency law.}
\end{figure}

The monkey model has a convoluted history.  The model can be thought of as
actually embedded within the combinatoric logic of Mandelbrot's early work providing an information-theoretic explanation of Zipf's law.  For example, in Mandelbrot \cite{Mandelbrot1954}     
there is an appendix entitled 
``Combinatorial Derivation of the Law $p_r=Pr^{-B}$'' (his notation) that effectively specifies a  monkey model, including a Markov version, but is not explicitly stated as such.  However, his derivation there is informal, and  the first completely rigorous analysis of the general monkey
model with independently combined letters was given only surprisingly recently by Conrad and Michenmacher \cite{C_and_M2004}.  Their analysis utilizes analytic number theory 
and is, as a result,  somewhat complicated - a point they note at the end of their article.  A simpler analysis using only elementary methods based on the Pascal pyramid has
now been given by Bochkarev and Lerner \cite{B_and_L2012}.  They 
have also analyzed the more general Markov problem \cite{B_and_L2014} and  hidden Markov models \cite{B_and_L2014b}.   
Edwards, Foxall and Perkins \cite{Edwards2012} have provided a directly relevant analysis in the context of scaling properties for paths on graphs 
explaining how the Markov variation can generate both a power law or a weaker scaling law, depending on the nature
of the transition matrix.  Perline \cite{Perline2015} also gives a simple approach to the demonstration of power law behavior.

A clarifying and important milestone in the history of this topic came from the cognitive psychologist,
G. A. Miller \cite{Miller1957,Miller1963}.     
  Miller used the simple model with equal letter probabilities and an independence assumption to show that this version generates a distribution of  word probabilities that conforms to an approximate  inverse power law such that the  probability $P_r$ of the $r^{\rm th}$ largest word probability is (in a sense that needs to be made precise) of the form $P_r\approx Cr^{-\beta}$, for $-\beta<-1$.  {\it In addition}, he made the very interesting observation that by using a keyboard of  $K=26$ letters and one space character, with the space character having a probability of .18 (about what is seen in empirical studies of English text), the value of $-\beta$ in this case is approximately -1.06, very close to the 
nearly -1 value  in samples of natural
 language text as illustrated in Figure 1.   In his simplified model, it turns out that $-\beta=-1 + \ln(1-s)/\ln K$, where $0<s<1$ is the space probability, so that
 $-\beta$ approaches -1 from below as $K$ increases.  (Natural logarithms are denoted $\ln(x)$ and logarithms with any other radix will be explicitly indicated, as in  $\log_K(x)$.  Note that Miller's  $\beta$ involves a  {\it ratio} of logarithms, so that it is invariant with respect to the radix used in the numerator and denominator.)
Consequently, Miller's model not only mimicks the form of Zipf's law for real languages, but with a large enough alphabet size it even mimicks the 
key parameter value.    {\it  In this article we give a broad generalization of Miller's observation to the case of unequal letter probabilities.}
In light of our present results given in Section 2, the application of sample spacings from a uniform distribution in \cite{Perline1996} to study the $-1$ behavior using an 
asymptotic regression line should now be viewed as just a step in a direction that ultimately led us to  the  Shao and Hahn limit used here.  

\section{Proof that \boldmath{$ -\beta$} Tends Towards   \boldmath{$-1$} from the Asymptotics of Log-Spacings}

For our analysis
we specify a keyboard with an alphabet of $K\ge 2$ letters $\{L_1\dots,L_K\}$ and a space character $S$.  The letter characters 
have non-zero probabilities $q_1,  q_2, \dots,  q_K$. The space character has probability $s$, so that $\sum_{i=1}^K q_i + s =1$.
A word is defined as any sequence of non-space letters terminating with a space.  
A word $W$ of exactly $n$ letters is a string such as  $W=L_{i_1}L_{i_2}\dots L_{i_n}S$ and    has a probability of the form $P(W)=P=q_{i_1}q_{i_2}\dots q_{i_n}s$ because letters are struck independently.  The space character with no preceding letter character will be considered a word of length zero.  
The rank ordered sequence of descending word probabilities  in the ensemble of all possible words is written
$P_1=s >  P_2\ge  \dots P_r  \ge \dots$ ($P_1=s$ is always the first and largest word probability.)  
 We  break ties for words with equal probabilities by alphabetical ordering, so that each word probability has a unique rank $r$. 

Conrad and Mitzenmacher \cite{C_and_M2004} give a carefully constructed definition of power law behavior in the monkey model as the situation where there exist two positive constants
$C_1$ and $C_2$ such that the inequality

\begin{equation}
C_1 r^{-\beta}\le P_r \le C_2 r^{-\beta}\ \ (\beta>0)
\end{equation}
holds for sufficiently large $r$.
Following the argument in Bochkarev and Lerner \cite{B_and_L2012},  $\beta$ in the monkey model turns out to be simply the solution to

\begin{equation}
q_1^{1/\beta} + q_2^{1/\beta} + \dots +q_K^{1/\beta}=1,
\end{equation}
with $\beta>1$.\footnote[1]{Bochkarev and Lerner \cite{B_and_L2012}  indicate  $\ln p(r) -\ln r/\gamma$ in  the  inequality
for their Theorem 1,  but  what was evidently intended is $\ln p(r) +\ln r/\gamma$.  Their $1/\gamma$ is equal to our $\beta$.} 
In the Miller model with equal letter probabilities, $q_1=\dots=q_K=(1-s)/K$, so   from equation (2), $\beta$ in this case is
found to be $1- \ln(1-s)/\ln K$.  In the Fibonacci example given by Conrad and Mitzenmacher \cite{C_and_M2004} and 
in Mitzenmacher \cite{Mitz2003}  
 they
use $K=2$ letters with probabilities $q_1$, $q_2=q_1^2$ and $q_1<(-1+\sqrt 5)/2$ so that $q_1+q_2<1$.   Then $\beta$ is  
the solution to  $q_1^{1/\beta}+q_1^{2/\beta}=1$, giving 

\begin{eqnarray}
\beta&=&\ln q_1 \over \ln ((-1 +\sqrt 5) /2).
\end{eqnarray}
To understand conditions leading to $-\beta\approx -1$,
we  define spacings through a random division of the unit interval and then state the Shao-Hahn limit law.  Let $X_1, X_2, \dots, X_{K-1}$ be a sample of $K-1$ i.i.d. random variables drawn from a distribution on $[0,1]$ with a bounded density function $h(x)\le M<\infty$.\footnote[2]{Shao and Hahn present the conditions for this limit in a more general way that reduces to this simpler statement when a density function exists.}
Write the order statistics of the sample as $X_{1:K-1}\ge X_{2:K-1}\dots \ge X_{K-1:K-1}$. The $K$ {\it spacings} $D_i$ are defined as the differences between the successive order statistics:
$D_1=1-X_{1:K-1}$, $D_i=X_{i-1:K-1}-X_{i:K-1}$ for $2\le i\le K-1$ and $D_K=X_{K-1:K-1}$.
We'll refer to this as a {\it generalized broken stick process}.
  By Shao and Hahn \cite{S&H}  Corollary 3.6, we have

\begin{equation}
{1\over K} \sum_{i=1}^K \ln(KD_i) \xrightarrow{a.s} -\int_0^1 h(x)\ln h(x) dx -\lambda,
\end{equation}  
as $K\to\infty$ and
where {\it a.s.} signifies {\it almost sure convergence}, $\lambda= .577\dots$ is the Euler constant and $-\int_0^1 h(x)\ln h(x) dx$ is the differential entropy of
$h(x)$.  
Clearly,

\begin{equation}
{1\over K} \sum_{i=1}^K \ln(KD_i) = \ln K +  {1\over K} \sum_{i=1}^K \ln D_i, 
\end{equation}  
so dividing through by $\ln K$ gives

\begin{eqnarray}
{\ln K  \over \ln K} +  {1\over K} {\sum_{i=1}^K \ln D_i \over \ln K} &\ &\nonumber \\
\xrightarrow{a.s} {-\int_0^1 h(x)\ln h(x) dx\over \ln K} &-&  {\lambda\over \ln K},
\end{eqnarray}
as $K\to\infty$.
The right side of the limit in (6) goes to 0 because ${-\int_0^1 h(x)\ln h(x) dx / \ln K}$ goes to 0 by the boundedness of the density $h(x)$.  Expressing
logarithms with a radix $=K$   then leads to the  limit

\begin{equation}
{\sum_{i=1}^K \log_K D_i \over K} \xrightarrow{a.s} -1\ \ {\rm as}\ K\to\infty.
\end{equation} 
Our {\it universality property} for $\beta$ will now follow almost immediately from this.
We use sample spacings to populate
the $K$ letter probabilities for the monkey keyboard.  
Since $s$ is the probability of the space character, define $q_i=(1-s)D_i$ $(1\le i \le K)$ so that $\sum_{i=1}^K q_i=1-s$.  
Let $\overline{m}_K=\sum_{i=1}^K \log_K q_i / K$.
Then from the limit (7), we  have 
\begin{eqnarray}
\overline{m}_K&=&\sum_{i=1}^K {\log_K q_i \over K}\nonumber \\
&=&{\log_K (1-s)} +{\sum_{i=1}^K \log_K D_i \over K}\xrightarrow{a.s} -1
\end{eqnarray}
as $K\to\infty$.
Since $\overline{m}_K\xrightarrow{a.s}-1$ and $-\beta<-1$, showing that $\overline{m}_K\le -\beta$ will prove that $-\beta\xrightarrow{a.s} -1$ as $K\to\infty$. 
  To see that this is the case, note:
\begin{eqnarray}
\sum_{i=1}^K {\log_K q_i ^{1/\beta} \over K}
&=& \log_K \Big[\bigl( q_1^{1/\beta} q_2^{1/\beta}\dots q_K^{1/\beta} \bigr)^{1/K}\Bigr]\\
&\le& \log_K \Big({ q_1^{1/\beta} +q_2^{1/\beta}\dots + q_K^{1/\beta}\over K}\Bigr)\\
&=&\log_K {1\over K}=-1,
\end{eqnarray}
where the inequality in (10) follows from the geometric-arithmetic mean inequality.  Therefore, from
\begin{eqnarray}
\sum_{i=1}^K {\log_K q_i ^{1/\beta} \over K}&=& {1\over \beta}\overline{m}_K \le -1,
\end{eqnarray}
we see that
 $\overline{m}_K\le -\beta$ and so $-\beta\xrightarrow{a.s.}-1$.  In the special case of Miller's model using all equal letter probabilities, 
$\overline{m}_K=-\beta$. 

The broad generality of the Shao-Hahn limit  leads to 
the near-universal characterization of the approximation $-\beta\approx -1$
as $K$ increases.  Figure 2 presents graphical results using simulations of sample spacings from several distributions to illustrate this phenomenon.
In Figure 2(a) we plot $\log_{10} P_r$ by $\log_{10} r$ for the Miller model
with exactly the parameter values he used, i.e., $K=26$ letters, a space probability $s=.18$ and equal letter probabilities $q_1=\dots=q_{26}=(1-.18)/26$.
The graph is based on the ranks $1\le r\le 475255=\sum_{i=0}^4 26^j$, corresponding to all words of length $\le 4$ non-space letters.
Figures 2(b) through 2(d) are based on generating $K=26$ letter probabilities from three different continuous distributions with bounded densities $h(x)$ on $[0,1]$ using  a generalized
broken stick method to obtain the sample  spacings.  Again $s=.18$ was used for the probability of the space character and the letter probabilities were populated with values
from the spacings with
$q_i=(1-s)D_i$, $1\le i \le 26$,
so that in each case, $\sum_{i=1}^{26} q_i=.82$, as in Miller's example.  For each graph in Figure 2(b) - 2(d), we generated 
 the largest 475255 word probability values in order to match the  Miller example of Fig 2(a).
 The three continuous distributions, all defined on $[0,1]$, are:
\begin{itemize}
\item[-] a unform distribution with density $h(x)=1$; 
\item[-] a beta $B(3,2)$ distribution with density 

\begin{equation}
h(x)={\Gamma(3+2)\over \Gamma(3)\Gamma(2)}x^{3-1}(1-x)^{2-1}; 
\end{equation}

\item[-] a triangular distribution with density
\begin{equation}
  h(x)=\left\{
  \begin{array}{@{}ll@{}}
    4x, & \text{if}\ 0\le 1/2 \\
    4(1-x) & \text{if}\ 1/2\le x\le 1.
  \end{array}\right.
\end{equation} 

\end{itemize}

\begin{figure}[h]
\centering
\includegraphics[width=4.in]{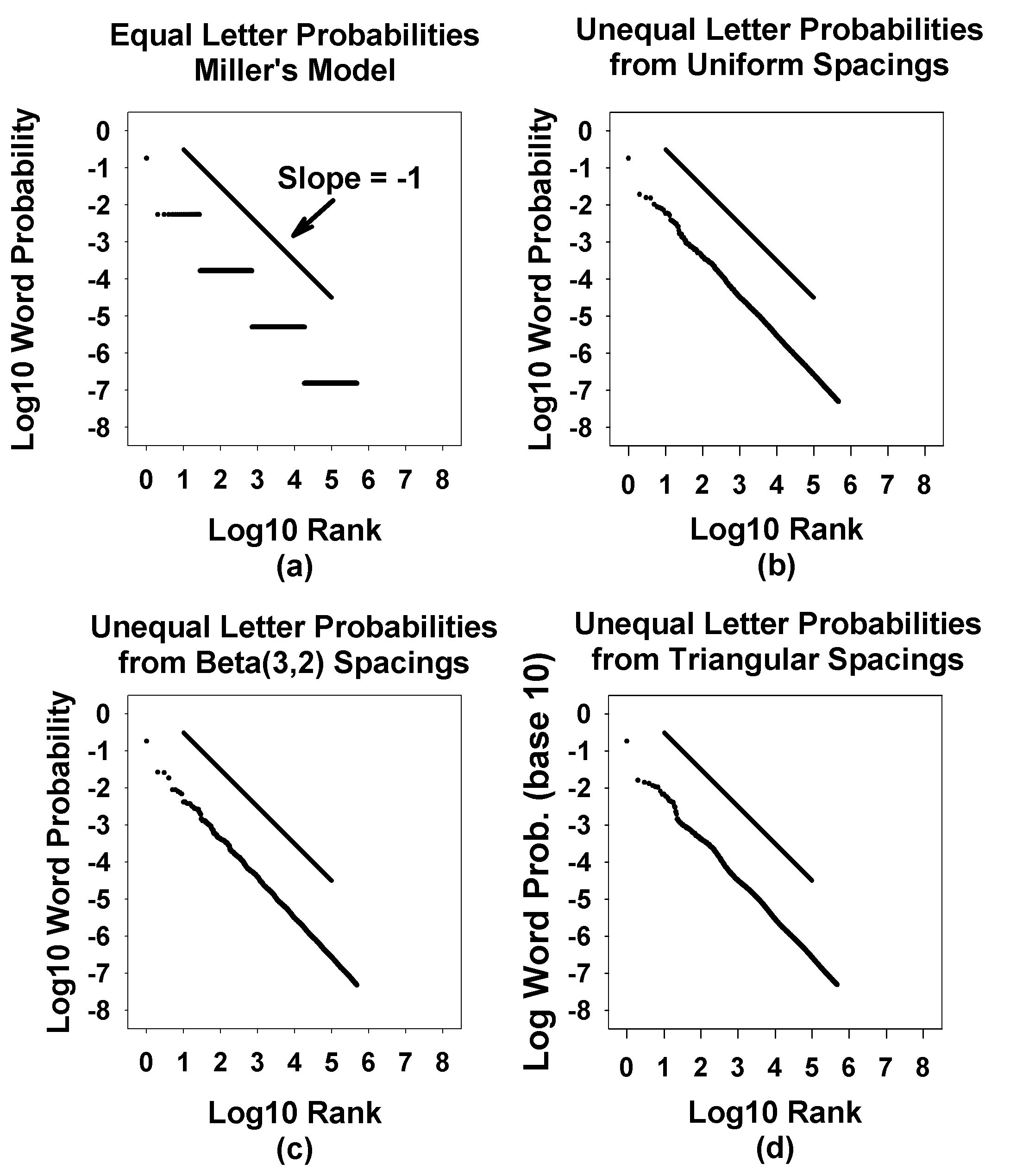}
\caption{\label{fig:PRL_Fig} Log-log plots of monkey word probabilities by rank showing  the asymptotic tendency towards a $-1$ exponent (slope on the log-log scales) using four different distributions to generate letter probabilities:
equal probabilities in 2(a) and a generalized broken stick process
in 2(b)-2(d).  $K=26$ letters were used in all cases and the largest $\sum_{i=0}^4 26^i=475255$ word probabilities are displayed. }
\end{figure}

The graphs in Figures 2(b)-2(d) illustrate our theoretical derivation,
but we also note that the very linear plots indicate what appears to be an almost ``immediate'' convergence
to power  law behavior exceeding what might be expected from the asymptotics of Conrad-Mitzenmacher and Bochkarov-Lerner.  

\section{Monkey Twitter:  Anscombe's CLT for the Model with a Finite Word Length Cutoff}

The discussion of the finite word length version of the monkey model in \cite{Perline1996} was incomplete because
it did not explicitly provide the normalizing constants for the application of Anscombe's CLT for a random number of
i.i.d. variables, as we do here.  We also want to explain more clearly the nature of the hybrid Zipf-lognormal distribution
that results from a finite length cutoff when letter probabilities are not identical.

The focus in Section 2
on the infinite ensemble of  word probabilities ${\bf P}_\infty$ {\it ranked}  to give $P_1>P_2\ge P_3\dots$ has actually
obscured  the significant {\it hierarchical structure} of the monkey model - what Mandelbrot \cite{Mandelbrot1983} (p. 345), called
a {\it lexicographic tree} - which only becomes evident when word length is considered.  To see this,  
write 
${\bf P}_{\rm len=n}$ for the multiset of all word probabilities for monkey words of exactly length $n$.
There are $K^n$ probabilities in ${\bf P}_{\rm len=n}$ having a total sum of $(q_1+q_2+\dots +q_K)^n  s=(1-s)^n s$ and an average value of
$((1-s)/K)^n s$, declining geometrically with $n$.
Now taking word length into account, Mandelbrot's lexicographic tree structure becomes clear with the simple case of $K=2$ letters:
the root node of the tree has a probability of $s$ (for the ``space'' word of length 0); it has two branches to the next level consisting of nodes for words of length 1 with probabilities
$q_1s$ and $q_2s$; these each branch out to the next level with four nodes having probability values $q_1^2 s, q_1q_2  s, q_2q_1s, q_2^2s$ and so on.   

The relevance of the CLT is seen first by
representing any probability $P\in{\bf P}_{\rm len=n}$  as a product of i.i.d.
random variables times the constant $s$:  $P=X_1X_2\dots X_n s$, where each $X_i$ takes on one of the letter probability values $q_1, q_2, \dots, q_K$
with probability $1/K$ (i.e., we use the natural counting measure to construct a probability space on ${\bf P}_{\rm len=n}$).
Let $\mu_1=\sum_{i=1}^K \ln q_i / K$ and $\sigma_1^2=\sum_{i=1}^K (\ln q_i - \mu_1)^2 / K$.  Assume that $\sigma_1^2>0$, i.e.,
the letter probabilities are not all equal.
Then $\mu_1 +\ln s$ is the mean and $\sigma_1^2$ is the variance of $\ln P$ for $P\in{\bf P}_{\rm len=1}$, and it follows that
$n\mu_1 +\ln s$ and $n\sigma^2_1$ are the respective mean and variance of $\ln P=S_n=\sum_{i=1}^n \ln X_i + \ln s$ for $P\in{\bf P}_{\rm len=n}$.    
Therefore, $(S_n - n\mu_1)/\sqrt{n\sigma_1^2}$ is asymptotically normally distributed $N(0,1)$ (the term $\ln s$ is rendered negligible in the asymptics)
 so that for sufficiently large $n$, $P$ itself will be approximately lognormal $LN(n\mu_1, n\sigma_1^2)$.
This approximate normality for the log-probabilities of words of fixed length is quite obvious and was noted in passing by Mandelbrot
\cite{Mandelbrot1961} (p. 210). However, he missed a {\it stronger} observation that the $N_n=\sum_{i=0}^n K^i$ word probabilities  for 
$P\in {\bf P}_{\rm len\le n}=\cup_{i=0}^n  {\bf P}_{\rm len=i}$  have a distribution that 
behaves in its central part away from the tails 
 very much like $P\in {\bf P}_{\rm len= n}$.   That is, a version of the CLT can be applied to words of length $\le n$, not just to words of length $n$, and
the two distributions, {\it in one sense}, are very close to each other.

To explain why this is so,  for $P\in{\bf P}_{\rm len\le  n}$ 
represent
$\ln P=S_{R_n}=\sum_{i=1}^{R_n} \ln X_i + \ln s$ as a sum of a {\it random number} of i.i.d. random variables (plus the constant $\ln s$), where $R_n$ is itself 
a random variable with a finite geometric distribution, ${\rm Prob}\{R_n=i\}=K^i / N_n, (1\le i \le n)$.  We will not repeat the
calculations of \cite{Perline1996} here, but  it is easily shown that $R_n/N_n\to 1$ in probability.
Because this limiting constant is 1,  Anscombe's generalization of the CLT \cite{Gut1987} (Theorem 3.1, p. 15) can be applied with the {\it identical}
normalizing constants $n\mu_1$ and $\sqrt{n\sigma^2_1}$  as used above with $\ln P$ for $P\in {\bf P}_{\rm len= n}$.  Consequently,  
it is {\it also} true that for $P\in {\bf P}_{\rm len\le n}$, the normalized sum
$(S_{R_n}-n\mu_1)/\sqrt {n\sigma_1^2}$ has an asymptotic $N(0,1)$ distribution.  In other words
the two random sums $S_n$ and $S_{R_n}$ behave so  similarly in their centers that 
 the word probabilities in  ${\bf P}_{\rm len= n}$ and
${\bf P}_{\rm len\le n}$ are {\it both} approximately $LN(n\mu_1,n\sigma_1^2)$.  However, it is essential to remark that
the two distributions have  very different upper tail behavior.  In this regard, 
Le Cam's \cite{LeCam1986} comment that French mathematicians use the term ``central'' referring to the CLT ``because it refers to the center of the distribution as opposed to its tails''
is particularly relevant.

The behavior of $P\in {\bf P}_{\rm len\le n}$ is illustrated in the graphs of Figure 3(a) and 3(b).  
The plot in Figure 3(a) was generated just as in Figure 2(b) except that a finite length cutoff of $n\le 4$ letters has been applied.
To make this clear, in Figure 2(b) the plot shows the top 475255 word probabilities in ${\bf P}_\infty$ generated using $K=26$ letters  derived from uniform spacings. 
In contrast, using the same letter probabilities, the plot in Figure 3(a) shows   the 475255 word probabilities generated with word length $\le 4$ letters, i.e.,
all the values of $P\in {\bf P}_{\rm len\le 4}$.
Word length in the case of {\it unequal} letter probabilities is certainly correlated, {\it but not perfectly}, with word probability:  
words of shorter length will, {\it on average}, have a higher probability than those of longer length, but except in Miller's degenerate case, there
will always be reversals where a longer word will have a higher probability than a shorter word.\footnote[3]{In natural languages, Zipf \cite{Zipf1949} underscored that ``the
length of a word tends to bear an inverse relationship to its relative frequency,'' which he called the {\it Law of Abbreviation.} 
 In many respects, this is the starting
point for his {\it Principle of Least Effort}. }
However, writing $P_{1:N_n}>P_{2:N_n} \dots \ge P_{N_n:N_n}$ for the ranked values of ${\bf P}_{\rm len\le n}$,
 it should be evident that for any given rank $r$,  $P_{r:N_n}\le P_r$  and that $P_{r:N_n}= P_r$ when $n$ is sufficiently large.  In short,
${\bf P}_{\rm len\le n}$ inherits its upper  tail power law behavior from ${\bf P}_{\infty}$, which is illustrated by
comparing the top part of the curve in Figure 3(a) with the corresponding part of Figure 2(b).  

Figure 3(b) uses the same data points from Figure 3(a) graphed  as a normal quantile plot, and its roughly  linear appearance for the logarithm of word probabilities for
 $P\in {\bf P}_{\rm len\le n}$ conforms to an approximate Guassian fit, although certainly the bending in the upper half of the distribution departs a bit from the linear trend of the lower half.  The fact that a distribution can have a power law tail and lognormal central part and {\it still look lognormal over essentially its   entire range}
may seem surprising.  The discussion in \cite{Perline2005} about {\it  power law mimicry} in the upper tail of lognormal distributions helps to explain
why this can happen.   We will call the distribution for  ${\bf P}_{\rm len\le n}$ a Zipf-lognormal distribution.

The monkey model with a fixed word length cutoff can prove useful as a motivating idea.  In the next section, we 
will discuss how models with the same branching tree structure have been proposed many times in the past, typically in settings where something like  
a finite word length is natural to consider.   For the moment, think of the social networking service, Twitter, which allows
members to exchange messages limited to at most 140 characters.  Define {\it Monkey Twitter} with a finite limit of $n+1$  characters.
For convenience,  in any implementation of a Monkey Twitter random experiment we will assume that monkeys would always fill up their allotted message space of $n+1$ characters.   Monkey words
still require a terminating space character, so it is possible for a monkey to type (at most) one ``non-word,'' which can vary in length from 1 to $n+1$ characters,
and will always be the last part of a message string.  Non-words are discarded, and the probabilities for the legitimate monkey words, varying in
length from 0 to $n$ non-space characters plus a terminating space character,  will correspond to the values in ${\bf P}_{\rm len\le n}$.

\begin{figure}[h]
\centering
\includegraphics[scale=.3]{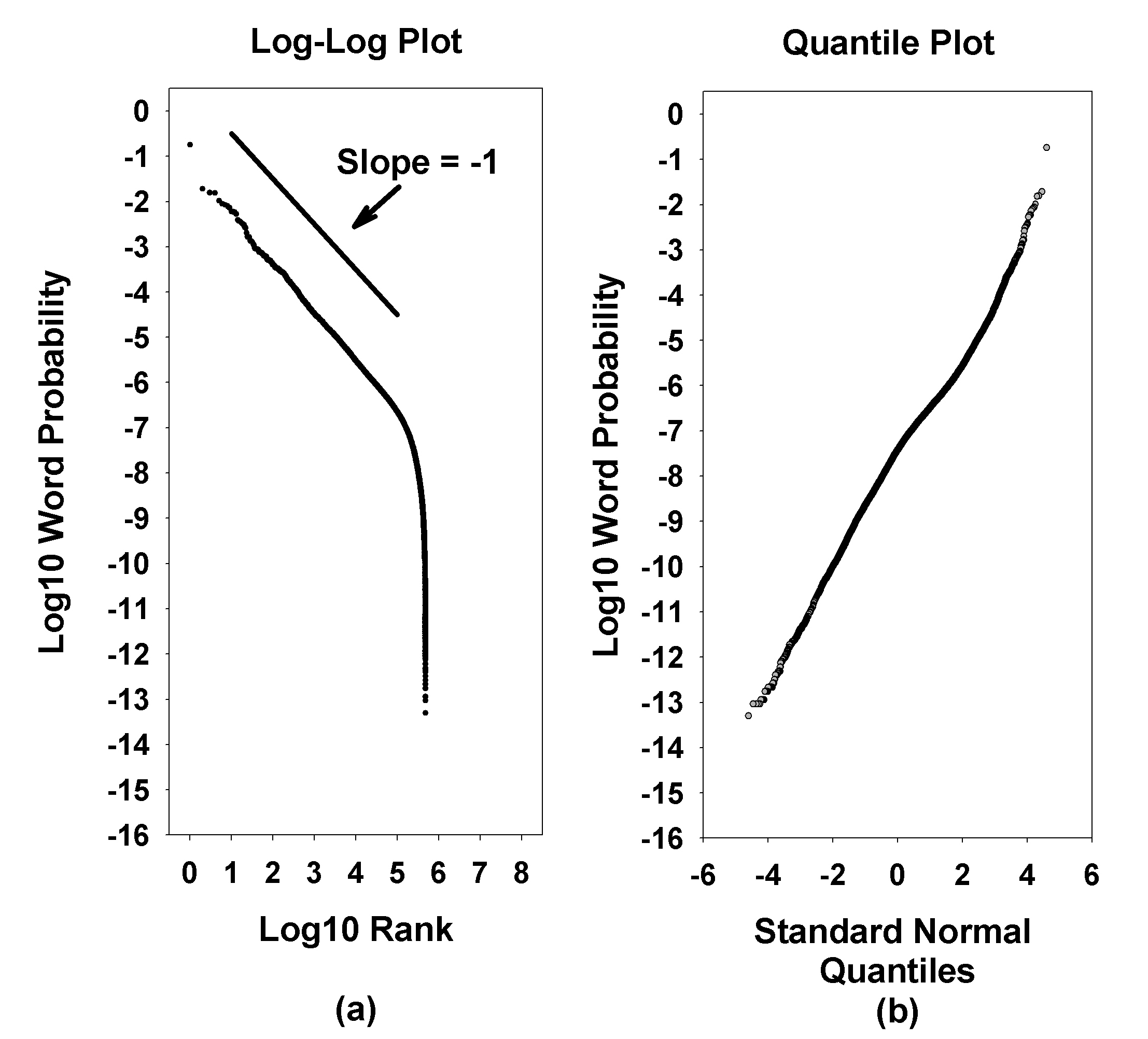}
\caption{ \label{fig:Rich3} Figure 3(a) is a log-log plot of monkey word probabilities by rank 
for all words of length $\le 4$ non-space characters using 
letter probabilities from uniform spacings. 
The linear upper tail coincides closely with the previous Figure 2(b), but  the power law clearly breaks down.  
Figure 3(b) shows the same word probabilities in a standard normal quantile plot.  Its rough linearity confirms an approximate  Gaussian fit over the whole distribution even 
though the upper tail is a power law as seen on the left in Figure 3(a).  We refer to this distribution as a Zipf-lognormal hybrid, and it has many
connections to other distributions discussed in the statistical literature. }
\end{figure}

\section{Connections to Other Work}

The subject of power law distributions is a vast topic of research reaching across diverse scientific disciplines \cite{Mitz2003, Perline2005, Clauset2009, Arnold2015}.   
Here we provide a brief sketch of how our results  connect to some other work and how they can be considered in a much more general light.  

The monkey model viewed in terms of its hierarchical tree structure is the starting point for understanding its general nature.
Though Miller did not explicitly
present his model as a branching tree, several researchers at almost the same time in the mid-to-late 1950s were highlighting this form using
the equivalent of equal letter probabilities
to motivate the occurrence of empirical power law distribution in other areas.  For example, Beckman \cite{Beckman1958} introduced this idea in connection with the power law of city populations
within a country (``Auerbach's law'' \cite{Auerbach1913}).  Using essentially the same logic as Miller, but in a completely different setting,  Beckman assumed that a given community will have $K$ satellite communities, each with a constant decreasing fraction of the population at a higher level.  That is, if $A$ is the population of the largest city and $0<p<1$, 
there will be $K$ nearby satellite communities with population $pA$, $K^2$ smaller communities nearby to these with population $p^2 A$, and so on. 
Mandelbrot \cite{Mandelbrot1997} (p. 226) had an apt expression for this pattern:  he described it
as ``compensation between two exponentials'' because the number of observations at each level increases geometrically while the mean value decreases geometrically.    
  Beckman then went on to note that
if instead of using a constant decreasing fraction $p$, one used a random variable $X$ on $(0,1)$, the population at the $n^{\rm th}$ level down
would be a random variable of the form $X_1 X_2 \cdots X_n A$, leading to an approximate lognormal distribution of populations at this level.  This corresponds to our discussion of monkey word probabilities in  ${\bf P}_{\rm len= n}$, although Beckman was not aware of the still stronger statement of lognormality for the probabilities of words of length $\le n$  
in ${\bf P}_{\rm len \le n}$ that we demonstrated  using Anscombe's  CLT.

Many more examples of a branching tree structure essentially equivalent to Miller's monkey model with equiprobable letters have been proposed over the years
to motivate the occurrence of a huge variety of empirical  power law distributions, including  such size distributions as lake areas, island areas (``Korcak's law''), river lengths, etc.  Indeed, Mandelbrot's
\cite{Mandelbrot1983} classic book on fractals is a rich source of these.  This equiprobability case is so simple to analyze
that it has been invoked over and over again   to illustrate how power law behavior can arise.  
  However, the more complicated case using unequal letter probabilities (or proportions) 
{\it and a finite word length cutoff} is what really underscores
the close analogy between the monkey Zipf-lognormal distribution and other mixture models exhibiting hybrid  power law tails and  approximately lognormal central ranges.

Empirical distributions of this form have appeared with great frequency.  In fact there is a long history
of researchers, including Pareto, first discovering what appears to be an empirical power law over an entire distribution because they start out by looking at only the largest observations (cities, corporations, islands, etc.).   However,  when they extend their measurements to the part of the distribution below the upper tail, {\it which is always more difficult to observe}, 
the power law behavior typically breaks down.  Several examples of this are chronicled in \cite{Perline2005}.   The power law for city sizes noted above is a perfect illustration.  Auerbach in 1913 \cite{Auerbach1913} looked at the
94 largest German cities from a 1910 census and showed a good power law fit, but because he did not include the thousands of  smaller communities, 
he may not have been aware that this relationship breaks down.  
Modern studies of the {\it full} distribution of communities, such as Eeckhout's  \cite{Eeckhout2004} analysis of the populations of 25,359 places from U.S. Census 2000, prove with high confidence that the {\it bulk} of the community populations fit an approximate lognormal distribution and that the power law behavior is confined to the upper tail.  

Montroll and Shlesinger \cite{Montroll1982, Montroll1983} explained how to generate a Pareto-lognormal hybrid  for income distributions
by using a hierarchical mixture model of lognormal distributions.
This is motivated from several angles, including the notion that higher classes amplify their income by organizing
 in such a way as to benefit from the efforts of others.  
Reed and Hughes \cite{Reed2002a}   have provided a far-ranging  framework
for understanding these hybrid distributions across the entire spectrum of disciplines
where they have been discovered:  ``physics, biology, geography, economics, insurance,
lexicography, internet ecology, etc.''   In \cite{Reed2002a} he and Hughes give
a condensed summary showing that ``if stochastic processes with exponential growth
in expectation are killed (or observed) randomly, the distribution of the killed or observed state 
exhibits power-law behavior in one or both tails.''  This work encompasses:  (1) geometric Brownian
motion (GBM); (2) discrete multiplicative processes; (3) homogeneous birth and death processes; and
(4) Galton-Watson branching processes.  In one of many variations on this theme, in his GBM model  \cite{Reed2004a} 
Reed specifies a stochastic process that uses lognormal distributions varying continuously in time with an exponential mixing distribution 
and shows 
that this generates  a ``Double Pareto-Lognormal Distribution.''   This has 
an asymptotic Pareto upper tail and a central part approximately lognormal;   in addition, it exhibits interesting asymptotic
behavior in the {it lower} tail where it is characterized by a direct (rather than an inverse) power law.  Reed has gone to great
effort to demonstrate the high quality of the fit of this distribution for all ranges of values (low, middle, high) across a wide variety of size distributions such as
particle and oil field sizes \cite{Reed2004a},  incomes \cite{Reed2004b},  internet file sizes and the sizes of biological genera \cite{Reed2002a}.

In today's nomenclature, the term {\it Zipf's law} has come to mean any power law distribution with an exponent close to the value $-1$, not just 
 the word frequency law.
(To be clear here, when we refer to power law distributions, we also mean the hybrids  with power law tails that we have been considering.)
The subset of power laws with this restricted exponent value  is surprisingly large and  includes the distribution of firm sizes \cite{Axtell2001},
city sizes \cite{Gabaix1999},   the famous Gutenberg-Richter law for earthquake magnitudes
\cite{Kagan1999} and many others.   

Gabaix \cite{Gabaix1999} pointed out that while it has long been known that stochastic growth processes could generate power laws,  ``these studies stopped short
of explaining why the Zipf exponent should be 1.''  To address this question, 
he has given a theoretical explanation of the genesis of a $-1$ exponent for  the Zipf law for city populations,
but in fact, his approach can  be regarded more generally. 
His key idea is a variation of Gibrat's \cite{Gibrat1931} {\it Law of Proportion Effect}: 
using a fixed number of  normalized city population quantities that sum to 1 and assuming growth processes (expressed as percentages)
randomly distributed as i.i.d. variables with a common mean and variance, he proves that a steady state will satisfy
Zipf's law. This proof requires
an additional strong assumption in order to reach a steady state, namely, a lower bound for the (normalized) population size.  
Gabaix goes on to discuss relaxed  versions of his model  and how it fits into the larger context of similar work done by others.

In their monograph on Zipf's law  Saichev et al \cite{Saichev2010} modify and extend the core Gabaix idea in numerous ways that render it more realistic.
For concreteness, they present their work using the terminology of financial markets and corporate asset values, but they are very clear on how
their models are relevant to a broad range of ``physical, biological, sociological and other applications.''   As with Reed,
GBM plays an important role in their investigations, which incorporate  birth and death processes, acquisitions and mergers,  the subtleties of finite-size effects
and other features.    Notably, they focus on the sets of conditions that lead to an approximate $-1$ exponent.  

In both Harr\"emoes and Tops{\o}e \cite{H_and_T2001} and Corominas-Murtra and Sol\'e \cite{C-M_and_Sole2010} $-1$ emerges as a consequence of
certain information theoretic ideas pertaining to the entropy function under growth assumptions.  In \cite{H_and_T2001} the focus is on language and the
 expansion of vocabulary size  while simultaneously maintaining finite entropy.  The perspective in \cite{C-M_and_Sole2010} is broader, but still 
based on the idea of systems with an expanding number of possible states ``characterized by special features on the behavior of the entropy.''  

There is another topic area of statistical research that  impacts on our work.  Natural language word frequency distributions have been characterized as
{\it Large Number of Rare Events} (LNRE) distributions because when collections of text are examined, no matter how big, there are always a 
large number of words that occur very infrequently in the sample \cite{Baayen2001}.
LNRE behavior indicates a tip-of-the-iceberg phenomenon resulting from a bias that captures the most common words, but
necessarily misses a vast quantity of rarely used words.
This is a classic and much studied problem encountered in species abundance surveys  in ecological research,
where a key question becomes estimating the size of the {\it zero abundance class} - the typically great number of species not observed in a sample
because of their rarity \cite{Bunge1993}.   

To get an idea of the significance of this issue, consider that Figure 3(a) graphs the {\it population} (or {\it parent} or {\it theoretical})
word probabilities for the Zipf-lognormal distribution, and not the {\it sample} frequencies as would be obtained from actually carrying out the Monkey Twitter experiment.
Simple visual inspection of Figure 3(a) indicates that the linear part of the graph (i.e., the power law tail) holds for about 5 logarithmic decades  or about
the first 100000 word probabilities out of the total of 475255 probabilities plotted there.  However, these 100000 probabilities   
comprise 97.4\% of the total probability mass of all 475255 values.   Intuitively, it should be clear that sampling from this Zipf-lognormal population distribution
will be dominated by the Zipf part, not the lognormal part.   Unless the sample size was astronomically large,  so that the large number of low probability words showed up,
the underlying structure of the parent distribution would not reveal itself.  To carry this idea still further, imagine the situation if the word length cutoff was on the
order of $n=140$ characters, such as with real Twitter.   No experiment could be run within any realistic time frame to ever hope to obtain a sample sufficiently large to uncover the true hybrid character of the population distribution -- the sample would always appear as a Zipf law, not a Zipf-lognormal law.  This kind of {\it visibility bias} has been a constant and recurring theme in all areas where Pareto-Zipf type distributions have been studied \cite{Perline1996, Perline2005}.  We believe its significance has been poorly appreciated in relation to this topic.

Finally, we will take a more bird's-eye view of matters and remark that our application of random spacings is very much in the spirit of the enormously fruitful study of {\it random
matrix theory} and a class of stochastic systems referred  to as the {\it KPZ universality class}.  Along these lines
 Borodin and Gorin \cite{Borodin2015}    have discussed a variety of probabilistic systems ``that can be analyzed by essentially algebraic methods,"
yet are applicable to a broad array of topics.
 Miller's demonstration of power law behavior and an approximate $-1$ exponent with increasing $K$ for the monkey model is in a similar vein.  
We regard this result as analogous to 
the Borodin and Gorin example of the De Moivre-Laplace proof of the CLT for a sequence of i.i.d. Bernoulli trials, which  
depends on having an explicit pre-limit distribution and then taking ``the limit of the resulting expression.''   (We want to add, however, that Miller's model is ultimately {\it too simple} to reveal the approximate lognormal behavior of the monkey word probabilities - for that, we needed to assume non-identical letter probabilities.)
Our point is to  make clear that the monkey model fits into a much larger conceptual scheme than appears at first glance.

\section{Conclusions}

Sample spacings provide  a natural way to populate the letter probabilities in the monkey model through a random division of the unit interval.
The Shao and Hahn asymptotic limit law 
for the logarithms of spacings then leads to the result that  the exponent in the power law of word frequencies will tend towards
a $-1$ value under broad conditions as the number of letters in the alphabet increases.   
The monkey model can also be viewed as a branching tree structure.  In that light, Anscombe's CLT 
reveals an underlying lognormal central part of the word frequency distribution
when a finite word length cutoff is imposed.  
The resulting hybrid Zipf-lognormal distribution has many connections to other work.
The visibility bias inherent in sampling from this distribution is similar to
what has been noted as a characteristic of many different empirical power laws.

\end{document}